\documentstyle{book}

\def\ru{\rule{-.5em}{0em}}

\begin{document}
\oddsidemargin 16.5mm
\evensidemargin 16.5mm
\thispagestyle{plain}

\noindent {\small\sc Univ. Beograd. Publ. Elektrotehn. Fak.}

\noindent {\scriptsize Ser. Mat. 7 (1996), 105--109.}

\vspace{5cc}
\begin{center}

{\Large\bf A NOTE ON HIGHER-ORDER\\
DIFFERENTIAL OPERATIONS\rule{0mm}{6mm}}

\vspace{1cc}
{\large \it Branko J. Male\v sevi\'c} 

\vspace{1cc}
\parbox{24cc}{{\scriptsize\bf In this paper we consider successive 
iterations of the first-order differential operations in space ${\bf R}^3.$}}
\end{center}

\vspace{1.5cc}
\begin{center}
{\bf 1. INTRODUCTION}
\end{center}

Let  ${\bf C}^{\infty} ({\bf R}^3)$ be the set of 
scalar functions $f=f(x_1 ,x_2 ,x_3 ):{\bf R}^3\mapsto {\bf R}$ which have  
the continuous partial derivatives of the arbitrary order on coordinates 
$x_i\; (i=1,2,3).$  Let $\vec {\bf C}^{\infty}({\bf R}^3)$ be the set vector 
functions $\vec f=\big( f_1(x_1,x_2,x_3),f_2(x_1,x_2,x_3),$\break 
$f_3(x_1,x_2,x_3)
\big):{\bf R}^3\mapsto {\bf R}^3$ which have the coordinately continuous
partial derivatives of the arbitrary order on coordinates $x_i\;\,(i=1,2,3).$
First-order differential operations of the vector analysis of the space
${\bf R}^3$ are defined on the following set of functions: 
$$ F=\left\{\, f:{\bf R}^3\mapsto {\bf R}\,| \, f\in {\bf C}^{\infty}({\bf R}^3)
\right\}\quad \mbox{and}\quad \vec F=\big\{ \vec f:{\bf R}^3\mapsto {\bf R}^3
\,|\, \vec f \in \vec {\bf C}^{\infty}({\bf R}^3)\big\}.$$

{\it First-order differential operations} of the vector analysis of the 
space ${\bf R}^3$ are defined as the following three linear operations 
[{\bf 1}], 
denoted here by $\nabla_1,\nabla_2$ and $\nabla_3$ for a convenience:  

{\small
\noindent (1)$\quad\displaystyle \mbox{grad}\,f=\nabla_1f= 
\frac{\partial f}{\partial x_1}\,\vec e_1+
\frac{\partial f}{\partial x_2}\,\vec e_2+
\frac{\partial f}{\partial x_3}\,\vec e_3: F\mapsto \vec F,$

\medskip
\noindent (2)$\quad\displaystyle \mbox{curl}\,\vec f=\nabla_2 \vec f= 
\left(\frac{\partial f_3}{\partial x_2}-
\frac{\partial f_2}{\partial x_3}\right)\vec e_1+
\left(\frac{\partial f_1}{\partial x_3}-
\frac{\partial f_3}{\partial x_1}\right)\vec e_2+
\left(\frac{\partial f_2}{\partial x_1}-
\frac{\partial f_1}{\partial x_2}\right)\vec e_3:\vec F\mapsto \vec F,$

\medskip
\noindent (3)$\quad\displaystyle \mbox{div}\,\vec f=\nabla_3 \vec f= 
\frac{\partial f_1}{\partial x_1}+
\frac{\partial f_2}{\partial x_2}+
\frac{\partial f_3}{\partial x_3}: \vec F\mapsto F.$}

\medskip
Let $\Omega =\{\nabla_1,\nabla_2,\nabla_3\}$ be the set of above defined 
operations and let $\Sigma =F \cup \vec F.$ Then the first-order differential 
operations can be considered as 
partial operations $\Sigma\mapsto \Sigma,$ i.e. as operations whose 
domain (and codomain) are subsets $F$ or $\vec F$ of $\Sigma.$
{\it Second and higher-order differential  operations} are 
then defined as products of operations in $\Omega$ in the sense of 
composition of operations. Some of these products might be {\it meaningful},
like $\nabla_3\circ \nabla_1,$ while the others are {\it meaningless}, 
like $\nabla_1 \circ \nabla_1.$ To all meaningless products for 
any argument we associate the value of nowhere defined function 
$\vartheta \;(\mbox{Dom}\,(\vartheta)=\emptyset$ and 
Ran\,$(\vartheta)=\emptyset).$  {\it Nowhere defined function}
$\vartheta (f_{\emptyset})$ is a concept from the 
recursive function theory [{\bf 2}]. We do not consider the function 
$\vartheta$ as the starting argument for calculating the value of the
higher-order differential operations. In that way we increase set 
$\Sigma$ into set $\Sigma =F\cup \vec F \cup \{\vartheta\}.$ 

All meaningful second-order differential operations are:

\medskip
\noindent $(4)\rule{5cc}{0cc}\Delta f= \mbox{div\,grad}\,f
=(\nabla_3 \circ \nabla_1)\,(f),$

\medskip
\noindent $(5)\rule{5cc}{0cc}\mbox{curl\,curl}\,\vec f
=(\nabla_2\circ \nabla_2)\,(\vec f),$

\medskip
\noindent $(6)\rule{5cc}{0cc}\mbox{grad\,div}\,\vec f
=(\nabla_1\circ \nabla_3)(\vec f),$

\medskip
\noindent $(7)\rule{5cc}{0cc}\mbox{div\,curl}\,\vec f
=(\nabla_3\circ \nabla_2)\,(\vec f)=0,$

\medskip
\noindent $(8)\rule{5cc}{0cc}\mbox{curl\,grad}\,f 
= (\nabla_2\circ \nabla_1)\,(f)=\vec 0,\quad 
f,\vec f\in \Sigma \setminus \{\vartheta \}.$

\medskip
In this paper we consider higher-order differential operations, 
search for meaningful ones and present some applications. 

\vspace{1cc}
\begin{center}
{\bf 2. HIGHER-ORDER DIFFERENTIAL OPERATIONS}
\end{center}

\medskip
\noindent {\bf Theorem 1.} {\it For arbitrary operations 
$\nabla_i,\nabla_j,\nabla_k \in \Omega \;\,(i,j,k\in \{1,2,3\})$ and 
argument $\xi\in \Sigma \setminus \{\vartheta \}$ the associative law holds:
$$ \nabla_i \circ (\nabla_j\circ \nabla_k) (\xi)=
(\nabla_i \circ \nabla_j)\circ \nabla_k (\xi).\leqno(9)$$}
\noindent {\bf Proof.} Choosing the $\nabla_i,\nabla_j,\nabla_k$
from $\Omega$ and argument $\xi$ from $\Sigma \setminus \{ \vartheta \},$
(9) appears in 54 possible cases. It is directly verified that whenever 
the left side of the equality is meaningless, the right side is also 
meaningless. Than, all meaningless products have the same value of the 
nowhere defined function $\vartheta,$ so that (9) is true in the following  
form: $\vartheta =\vartheta.$ Also, whenever the left side of equality is
meaningful, the right side is also meaningful. Then, according to the 
associative law of the meaningful functions, we conclude that (9) is true. 

From Theorem 1 it follows (by induction) that the generalized
associative law also holds, so we may write the product 
$\nabla_{i_1}\circ \nabla_{i_2}\circ \cdots \circ \nabla_{i_n}$   
without brackets $(i_j\in\{1,2,3\}:\; j=1,2,...,n).$ 

For higher-order differential operations, given as meaningful 
products, we say that they are the {\it trivial products} if they are
trivially anullated, i.e. if they are identically the same as the anullating 
functions $0,\,\vec 0$ from $\Sigma.$ Otherwise, we refer to the higher-order 
differential operations, given as meaningful products, as {\it nontrivial 
products} (if they are nontrivially anullated). 

Next, we prove the statement: 

\noindent {\bf Theorem 2.} {\it Higher-order differential operations
appear as nontrivial products in the following three forms}: 
\begin{eqnarray*}
\mbox{(grad)\,div\,\ldots\,grad\,div\,grad}\, f 
\ru & = & \ru (\nabla_1\circ) \nabla_3 \circ \cdots \circ \nabla_1\circ\nabla_3
\circ \nabla_1 f,\\
\mbox{curl\,curl\,\ldots\,curl\,curl\,curl}\,\vec f 
\ru & = & \ru \nabla_2\circ \nabla_2 \circ \cdots \circ \nabla_2\circ\nabla_2
\circ \nabla_2 \vec f,\\
\mbox{(div)\,grad\,\ldots\,div\,grad\,div}\,\vec f
\ru & = & \ru (\nabla_3\circ) \nabla_1 \circ \cdots \circ \nabla_3\circ\nabla_1
\circ \nabla_3 \vec f,
\end{eqnarray*}
{\it for arbitrary functions $f,\,\vec f\in \Sigma\setminus \{ \vartheta \},$
where terms in brackets are included for odd number of terms and are left 
out otherwise. All other meaningful operations are identically zero 
in their domain.}

\noindent {\bf Proof.} Meaningful third-order differential operations appear
in the form of eight compositions as follows: 

\medskip
\noindent$(10)\rule{5cc}{0cc} \mbox{grad\,div\,grad}\, f
=\nabla_1\circ \nabla_3 \circ \nabla_1 f,$

\medskip
\noindent$(11)\rule{5cc}{0cc} \mbox{curl\,curl\,curl}\,\vec f
=\nabla_2\circ \nabla_2 \circ \nabla_2 \vec f,$

\medskip
\noindent$(12)\rule{5cc}{0cc} \mbox{div\,grad\,div}\,\vec f
=\nabla_3\circ \nabla_1 \circ \nabla_3 \vec f,$

\medskip
\noindent$(13)\rule{5cc}{0cc} \mbox{div\,curl\,curl}\,\vec f
=\nabla_3\circ \nabla_2 \circ \nabla_2 \vec f=0,$

\medskip
\noindent$(14)\rule{5cc}{0cc} \mbox{div\,curl\,grad}\, f
=\nabla_3\circ \nabla_2 \circ \nabla_1 f=0,$

\medskip
\noindent$(15)\rule{5cc}{0cc} \mbox{curl\,curl\,grad}\, f
=\nabla_2\circ \nabla_2 \circ \nabla_1 f=\vec 0,$

\medskip
\noindent$(16)\rule{5cc}{0cc} \mbox{curl\,grad\,div}\,\vec f
=\nabla_2\circ \nabla_1 \circ \nabla_3 \vec f=\vec 0,$

\medskip
\noindent$(17)\rule{5cc}{0cc} \mbox{grad\,div\,curl}\,\vec f
=\nabla_1\circ \nabla_3 \circ \nabla_2 \vec f=\vec 0, \quad
f,\,\vec f\in \Sigma \setminus \{\vartheta \}.$

\medskip
\noindent Anullations of the operations (13)--(17) follow directly from the 
anullations (4)--(5). The statement follows directly from the principle 
of mathematical induction by means of using the general associative law 
and formulas (10)--(17). 

For a given sequence of operations $\nabla_{i_1},\nabla_{i_2},
\ldots ,\nabla_{i_n}$ from the set $\Omega$ of functions, let define 
the concept of the {\it collection of functions} as a subset of functions
$\Theta \subseteq \Sigma \setminus \{ \vartheta \}$ such that all 
functions $\xi$ from $\Theta$ anullate the nontrivial product 
$\nabla_{i_1}\circ \nabla_{i_2}\circ \cdots \circ \nabla_{i_n}
\,(\xi).$

Let us form some collections. Scalar functions $f$ from $\Sigma,$
such that $\Delta^n f=0$ is true, define {\it harmonic collection}
${\bf H}_n$ of order $n,$ as the form of the polyharmonic functions.
Let us notice that in the case of two dimensions there is a general
form of polyharmonic functions $f$ as a solution of the equation 
$\Delta^n f=0,$ [{\bf 3}]. Vector functions $\vec f$ from $\Sigma,$ such that 
curl$^n \vec f = \vec 0$ is true, define {\it curling collection}
${\bf C}_n$ of order $n.$

We can remark that besides the total scalar operation
$\Delta:F\mapsto F$ (partial scalar operation $\Delta : \Sigma 
\mapsto \Sigma $) we can also consider the total vector operation 
$\vec \Delta :\vec F \mapsto \vec F$ (partial vector operation
$\vec \Delta :\Sigma \mapsto \Sigma$) defined by: 
$$ \vec \Delta \vec f=(\Delta f_1,\Delta f_2,\Delta f_3)
=\Delta f_1\cdot \vec e_1 +\Delta f_2\cdot \vec e_2
+ \Delta f_3\cdot \vec e_3. \leqno(18)$$
Let set \bf $\vec{\bf H}_{n}$ \rm be the sign for the vector functions 
$\vec{f}$ from $\Sigma$ such that $\vec{\Delta}^{n}(\vec{f})=\vec{0}$, 
where $\vec{\Delta}^{n}$ is iteration of order n of the vector operation 
$\vec{\Delta}$ given by  (18). \it The set of vector harmonic functions 
\bf $\vec{\bf H}_{n}$ \rm of order n, which is defined in such a way, is not 
in the list of collections which appear in the previous theorem because 
it is not obtained through the compositions of operations (1)--(3). For 
the set \bf $\vec{\bf H}_{n}$ \rm we shall keep the term {\it collection.}

Let us notice that for scalar polyharmonic collections, vector
polyharmonic collections and curling collections, related to the index-order, 
the following inclusions hold: 
$$ {\bf H}\subset {\bf H}_2 \subset \cdots \subset {\bf H}_{n-1}
\subset {\bf H}_n\subset \cdots,\leqno(19)$$
\vspace{-6mm}
$$ \vec{\bf H}\subset \vec{\bf H}_2 \subset \cdots \subset \vec{\bf H}_{n-1}
\subset \vec{\bf H}_n\subset \cdots,\leqno(20)$$
\vspace{-6mm}
$$ {\bf C}\subset {\bf C}_2 \subset \cdots \subset {\bf C}_{n-1}
\subset {\bf C}_n
\subset \cdots\;.\leqno(21)$$

Let emphasize that all previous considerations can be transformed in 
three-dimensional orthogonal curvilinear coordinate system by introducing 
of corresponding presumptions for functions from the sets $F,\,\vec F$
and {\sc Lam\'e}'s coefficients. 

Finally, let state a few examples where scalar and vector polyharmonic
collections appear. 

\medskip
\noindent 
{\bf Example 1.} {\it All meaningful products of third-and-higher-order
differential operations for vector functions $\vec f\in \vec {\bf H}$
and scalar functions $f\in {\bf H}$ are anullated.} 

\smallskip
\noindent 
For vector functions $\vec f\in \vec{\bf H}$ the following equation holds:
$$ \mbox{curl\,curl}\,\vec f = \mbox{grad\,div}\,\vec f.\leqno(22)$$ 
Hence, for $f\in {\bf H}$ and $\vec f\in \vec{\bf H},$ on the basis 
of formulas (22) and (10)--(17) the following is true: 
\begin{eqnarray*}
& & \mbox{grad\,div\,grad}\,f = \mbox{grad}\,(\Delta f)=\vec 0,\\                     
& & \mbox{curl\,curl\,curl}\,\vec f =  \mbox{curl}\,(\mbox{grad\,div)}\,
\vec f = \vec 0,\\
& & \mbox{div\,grad\,div}\,\vec f =  \mbox{div\,(curl\,curl)}\,
\vec f= 0.                 
\end{eqnarray*}
Thus, all eight meaningful products of third-order differential operations 
are anullated, so that the statement is true. 

\medskip
\noindent 
{\bf Example 2.} {\it If $f\in {\bf H}_{n-1},$ then
$x\cdot f\in {\bf H}_n,\;\;n\ge 2.$}

\smallskip
\noindent 
Let us notice that if $f\in F,$ then $x\cdot f\in F.$ For an arbitrary
scalar function $f\in F$ the following equation is directly verified: 
$$ \Delta (x\cdot f)=2\partial f/\partial x+x\cdot \Delta (f).$$
Inductive generalization is the following equation:
$$ \Delta^n(x\cdot f)=2n\cdot \partial \big( \Delta^{n-1}(f)\big)/\partial x
+x\cdot \Delta^n(f).$$
Thus, for $(n-1)$-harmonic function $f\in {\bf H}_{n-1}$ the conclusion 
$x\cdot f\in {\bf H}_n$ is true.

\break

\medskip
\noindent 
{\bf Example 3.} {\it If $f\in {\bf H}_{n-1},$ then 
$(x^2+y^2+z^2)\cdot f\in {\bf H}_n,\;\,n\ge 2.$}

\smallskip
\noindent 
Let us notice that if $f\in F,$ then $(x^2+y^2+z^2)\cdot f\in F.$ For the
arbitrary scalar function $f\in F$ the following equations are directly 
verified: 
\begin{eqnarray*}
\Delta (x^2\cdot f)\ru & = & \ru 2\cdot f+4x\cdot \partial f/\partial x
+x^2\cdot \Delta(f),\\
\Delta^2(x^2\cdot f)\ru & = & \ru 8\cdot \partial^2 f/\partial x^2
+8x\cdot \partial \big( \Delta (f)\big)/\partial x
+4\cdot \Delta (f) +x^2\cdot \Delta^2 (f).
\end{eqnarray*}
Inductive generalization is the equation as follows: 
\begin{eqnarray*}
& & \Delta^n(x^2\cdot f)=4n(n-1)\cdot \partial^2 \big(\Delta^{n-2}(f)\big)/
\partial x^2\\
& & \rule{5cc}{0cc} + 4nx\cdot \partial \big( \Delta^{n-1}(f)\big)/\partial x
+ 2n\cdot \Delta^{n-1}(f) +x^2\cdot \Delta^n (f).
\end{eqnarray*}
Thus, if $f\in {\bf H}_{n-1},$ then $(x^2+y^2+z^2)\cdot f\in {\bf H}_n.$ 

\bigskip 
\noindent 
Two previous examples are the generalizations of 
the corresponding problems contained in [{\bf 4}].

\bigskip
\noindent 
{\bf Acknowledgement.} I wish to express my gratitude to Professors 
{\sc M. Merkle, I. Lazarevi\'c} and {\sc D. To\v si\'c} who examined 
the first version of paper and gave me their suggestions and some very 
useful remarks. 

\vspace{2cc}
\begin{center}{\small\bf REFERENCES}
\end{center}

\vspace{0.8cc}
\newcounter{ref}
\begin{list}{\small \arabic{ref}.}{\usecounter{ref} \leftmargin 4mm 
\itemsep -1mm}

\item {\small {\sc M. L. Krasnov, A. I. Kiselev, G. I. Makarenko:} 
{\it Vector Analysis.} Moscow 1981.}

\item {\small {\sc N. Cutland:} {\it Computability.}
Cambridge University Press, London 1980.} 

\item {\small {\sc D. S. Mitrinovi\'c, J. D. Ke\v cki\'c:} 
{\it Jedna\v cine matemati\v cke fizike.} Beograd 1985.}

\item {\small {\sc D. S. Mitrinovi\'c}, in association with 
{\sc P. M. Vasi\'c}: {\it Diferencijalne jedna\v cine}, 
Novi zbornik problema 4. Beograd 1986.} 

\item {\small {\sc M. J. Crowe:} {\it A History of Vector Analysis.} 
University of Notre Dame Press, London 1967.} 
\end{list}

\vspace{1 cc}

{\small
\noindent Faculty of Electrical Engineering,  \hfill 
                              (Received May 6, 1996)\break
\noindent University of Belgrade,             \hfill\break
\noindent P.O.B 816, 11001 Belgrade,          \hfill\break
\noindent Yugoslavia                          \hfill\break
\noindent {\bf malesevic@kiklop.etf.bg.ac.yu} \hfill\break}

\newpage

\end{document}